\documentclass[preprints,article,accept,moreauthors,LaTeX, dvi2pdf]{mdpi} 

\usepackage{bm}
\usepackage{mathtools}
\usepackage{siunitx}
\usepackage{algorithm}
\usepackage{algorithmic}
\usepackage{siunitx}
\usepackage{tabularx}
\usepackage{subcaption}
\usepackage{setspace}

\newtheorem{thm}{Theorem}
\newcommand{\D}{\mathrm{d}}

\firstpage{1} 
\pubvolume{xx}
\issuenum{x}
\articlenumber{xx}
\pubyear{2020}
\copyrightyear{2020}
\history{}

\Title{Mathematical derivation for Vora-Value based filter design method: Gradient and Hessian }

\Author{Yuteng Zhu $^{1}$* and Graham D. Finlayson $^{1}$}
\AuthorNames{Yuteng Zhu and Graham D. Finlayson}
\address{%
$^{1}$ \quad School of Computing Sciences, University of East Anglia}
\corres{Correspondence: yuteng.zhu@uea.ac.uk}
\firstnote{Current address: Norwich Research Park, Norwich, NR4 7TJ, U.K.} 

\abstract{In this paper, we present the detailed mathematical derivation of the gradient and Hessian matrix for the Vora-Value based colorimetric filter optimization. We make a full recapitulation of the steps involved in differentiating the objective function and reveal the positive-definite Hessian matrix  when a positive regularizer is applied. This paper serves as a supplementary material for our paper in the colorimetric filter design theory. }

\begin{document}

\section{Preliminary}
Let $\mathbf{Q} = \mathbf{[r,g,b]}$ and  $\mathbf{X} = \mathbf{[x,y,z]}$ denote respectively the spectral sensitivities of the camera and the CIE XYZ color matching functions (CMFs) of the human visual sensors. The columns of matrices $\mathbf{Q}$ and  $\mathbf{X}$ represent the spectral sensitivity for each sensor channel and the rows represent the sensor responses at a sampled wavelength. Both matrices are in the size of $n \times 3$, where $n$ is the number of sampling wavelengths across the visible spectrum.

\subsection{Notation}
We will use the following notation for the gradient and Hessian matrix (with respect to the $n$-dimensional filter vector $\mathbf{f} = [\mathbf{f}_1, \mathbf{f}_2,\cdots, \mathbf{f}_{n}]^T$):
\begin{equation}
\nabla \nu (\mathbf{f}) =  
\begin{bmatrix}
    \frac{\partial \nu}{\partial \mathbf{f}_1}\\
    \frac{\partial \nu}{\partial \mathbf{f}_2}\\
    \vdots \\
    \frac{\partial \nu}{\partial \mathbf{f}_{n}}\\
\end{bmatrix}
\hspace{0.25in} \text{and} \hspace{0.25in}
   \mathbf{ H} = \nabla ^2 \nu  (\mathbf{f}) = 
    \begin{bmatrix}
        \frac{\partial^2 \nu}{\partial \mathbf{f}_1^2} & \frac{\partial^2 \nu}{\partial \mathbf{f}_1 \mathbf{f}_2} & \cdots & \frac{\partial^2 \nu}{\partial \mathbf{f}_1 \mathbf{f}_{n}}\\
        \frac{\partial^2 \nu}{\partial \mathbf{f}_2 \mathbf{f}_1} & \frac{\partial^2 \nu}{\partial \mathbf{f}_2^2} & \cdots & \frac{\partial^2 \nu}{\partial \mathbf{f}_2 \mathbf{f}_{n}} \\
        \vdots & \vdots & \ddots & \vdots \\
        \frac{\partial^2 \nu}{\partial \mathbf{f}_{n} \mathbf{f}_1} & \frac{\partial^2 \nu}{\partial \mathbf{f}_{n}\mathbf{f}_2} & \cdots & \frac{\partial^2 \nu}{\partial \mathbf{f}_{n}^2}
\end{bmatrix}
\end{equation}
or equivalently, using the indices, we can respectively express as $(\nabla \nu)_i = \frac{\partial \nu}{\partial \mathbf{f}_i}$ and $\mathbf H_{i,j} = \frac{\partial ^2 \nu}{\partial \mathbf{f}_i \mathbf{f}_j}$~\cite{bibmatrix}.
From the definitions, we know that $\nabla \nu$ and $\nabla ^2 \nu$ are respectively in the size of $n \times 1$ and $n \times n$.

\subsection{Vora-Value}
Given a camera sensor set $\mathbf{Q}$ and the trichromatic human visual sensors $\mathbf{X}$, the Vora-Value is defined~\cite{bibVV} as
\begin{equation}
\begin{split}
     \nu (\mathbf{ Q,X}) & = \frac{1}{3} tr(\mathbf{{Q(Q^T Q)^{-1} Q^T X(X^T X)^{-1}X^T}})
\end{split}
\label{eq:vv_def}
\end{equation}
where the superscripts $^T$ and $^{-1}$ denote respectively the matrix transpose and inverse and $tr()$ returns the  sum of the elements along the diagonal of a matrix. The Vora-Value is often used to measure  how similarly a camera samples the spectral signals compared to the human visual system. It returns a number in the range [0,1] where 1 means the camera is fully colorimetric such that RGBs are precisely a linear transform from XYZ tristimulus values. A higher Vora-Value indicates a better fit between two sensor systems.

\subsection{Projector Matrix}
The projector of a matrix - such as $\mathbf{Q}$ - is defined as
\begin{equation}
    \mathrm{P}\{\mathbf{Q}\}=\mathbf{Q(Q^T Q)^{-1}Q^T}.
    \label{eq:proj}
\end{equation}
When we return to Eq.~(\ref{eq:vv_def}), it can also be written in a more compact representation as 
\begin{equation}
\begin{split}
     \nu (\mathbf{Q,X})&=\frac{1}{3}tr(\mathrm{P}\{\mathbf{Q}\}\mathrm{P}\{\mathbf{X}\})
\end{split}
\label{eq:vv_def_proj}
\end{equation}
where $\mathrm{P}\{\mathbf{Q}\}$ and $\mathrm{P}\{\mathbf{X}\}$ denote the projection matrices respectively of the camera spectral sensitivities $\mathbf{Q}$ and the human visual responses $\mathbf{X}$.

\subsection{Orthonormal Basis}

In a $n$-dimensional vector space, $n$ linearly independent vectors forms a set. We call such a set as {\em basis set}. Every vector in the $n$-dimensional vector space can be expressed as a linear combination of the basis vectors. There are infinite bases for a vector space and, by definition, they are all linear transform apart. 

Let $\mathbf{V = [v_1,v_2, v_3]}$ denote a special linear combination of $\mathbf{X} = \mathbf{[x,y,z]}$ as
\begin{equation}
    \mathbf{V}=\mathbf{XT}
\label{eq:linearMapping}
\end{equation}
where $\mathbf{T}$ is the (full rank) linear mapping matrix which makes $\mathbf{V}$ orthonormal. An orthonormal matrix has columns that are unit vectors and also perpendicular to each other. Mathematically, we write $\mathbf{V}^T \mathbf{V} = \mathbf{I}_{3}$ ($\mathbf{I}_{3}$ is the ${3\times 3}$ identity matrix). The orthonormal matrix $\mathbf{V}$ can be obtained by many methods, e.g.\ the Gram-Schmidt process~\cite{bibGolub}.
 
By simple substitution into the matrix projector in Eq.  (\ref{eq:proj}), we can express projector matrices in a simpler algebraic form:
\begin{equation}
\mathrm{P}\{\mathbf{X}\} = \mathrm{P}\{\mathbf{V}\} =\mathbf{VV^T}
\label{eq:proj_X}
\end{equation}
and then substituting into Eq.~(\ref{eq:vv_def}), we can simplify the Vora-Value as
\begin{equation}
\begin{split}
     \nu (\mathbf{Q,X}) = \nu (\mathbf{Q,V}) 
     =\frac{1}{3}tr(\mathbf{Q(Q^T Q)^{-1} Q^T VV^T} ).
\end{split}
\label{eq:vv_def_simple}
\end{equation}

\subsection{Filter-modified Vora-Value}
Previously, we proposed to design a color filter which, when placed in front of a camera, can make the new effective camera more colorimetric~\cite{bibLutherFilter}. When a color filter is placed in front of a camera, it alters the spectral sensitivities. The effect of placing a color filter to a camera can be modeled as the multiplication of the filter spectral transmittance to the camera spectral sensitivities. Given $\mathbf f$ an $n$-dimensional filter vector (with $\mathbf{f}_i > 0$) and camera spectral sensitivity matrix $\mathbf Q$, the new effective sensitivity responses after filtering can be represented as $\mathbf{diag(f)Q}$. To ease the notation, we use $\mathbf{F = diag(f)}$ and rewrite as $\mathbf{FQ}$. 

Thus the filter-modified Vora-Value for the effective `filter+camera' system (using the orthonormal basis of the XYZ CMFs) can be written as
\begin{equation}
\begin{split}
 \nu (\mathbf{FQ, X}) &= \frac{1}{3} tr(\mathbf{FQ(Q^T F^2Q)^{-1}Q^T F \, VV^T}).
\end{split}
\label{eq:vv_fcam}
\end{equation}
Or equivalently, in the simpler representation (using projector matrices), we write $ \nu (\mathbf{FQ, X})= \frac{1}{3} trace(\mathrm{P}\{\mathbf{ FQ}\}\mathrm{P}\{\mathbf{X}\})  $.

\section{Derivation of Gradient}
In this section, we will derive the gradient, in terms of the filter vector, of the filter-modified Vora-Value as given in Eq.~(\ref{eq:vv_fcam}). 
\begin{thm}\label{theorem:gradient}
$\nabla \nu ({\mathbf f}) =  \frac{\partial \nu(\mathbf{F}) }{\partial \mathbf{f}} = \frac{2}{3} \, ediag \Big( \mathbf{F}^{-1}\mathrm{P}\{ \mathbf{FQ}\} \, \mathrm{P}\{\mathbf{X}\} (\mathbf{I} -\mathrm{P}\{ \mathbf{FQ}\}) \Big)$
\end{thm}

\noindent
\begin{proof}
The following rules of matrix calculus are used to obtain the required differentials:
\begin{equation}
    \begin{split}
    \D\, tr(\mathbf U) &= tr(\D \mathbf U) \\
    \D(\mathbf{UV}) &= \mathbf{U \, \D \mathbf{V} +} \D \mathbf{U} \, \mathbf{V} \\
    \D (\mathbf{AU}) &= \mathbf{A} \,\D \mathbf{U} \\
    \D \mathbf{\mathbf{U}^{-1}} & = -\mathbf{U^{-1}} (\D \mathbf{U}) \mathbf{U^{-1}} 
    \end{split}
    \label{eq:calculus_rules}
\end{equation}

\noindent
Using the above rules, we have
\begin{equation}
\begin{split}
    \D \nu(\mathbf{F}) = & \, \frac{1}{3}\, tr \Big( \D \mathbf{F \, Q(Q^T F^2Q)^{-1}Q^T F \, VV^T -}\\
     &2  \mathbf{FQ (Q^T F^2Q)^{-1} Q^T F \,\D F \, Q (Q^T F^2Q)^{-1} Q^T F \, VV^T  + }\\
     &\mathbf{FQ (Q^T F^2Q)^{-1}Q^T \, \D F \, VV^T} \Big).
\end{split}    
\label{eq:d_nu}
\end{equation}
\noindent
Using the acyclic property of trace that $tr(\mathbf{ABC}) = tr(\mathbf{BCA}) = tr(\mathbf{CAB})$, we can move the $\D \mathbf F$ in each of the term to the end of the formulation. We also use the projector representation of $\mathrm{P}\{\mathbf{X}\} = VV^T$ to make it more compact as
\begin{equation}
\begin{split}
    \D \nu(\mathbf{F}) = & \, \frac{1}{3}\, tr \Big( \mathbf{Q(Q^T F^2Q)^{-1}Q^T F \, \mathrm{P}\{\mathbf{X}\} \D F}  - \\
& 2 \mathbf{Q (Q^T F^2Q)^{-1} Q^T F \,\mathrm{P}\{\mathbf{X}\} \, FQ (Q^T F^2Q)^{-1}Q^T F \,\D F} + \\
& \mathbf{\mathrm{P}\{\mathbf{X}\} FQ (Q^T F^2Q)^{-1}Q^T 
\D F} \Big).
\end{split}
\end{equation}

\noindent
The matrix $\mathbf{F}$ is diagonal by definition, we can make the derivative as
\begin{equation}
\begin{split}
\frac{\partial \nu(\mathbf{F}) }{\partial \mathbf{F}_{ii}}  = & \frac{1}{3} \Big[ \mathbf{Q(Q^T F^2Q)^{-1}Q^T F \, \mathrm{P}\{\mathbf{X}\}}  - \\
& 2 \mathbf{Q (Q^T F^2Q)^{-1} Q^T F \,\mathrm{P}\{\mathbf{X}\} \, FQ (Q^T F^2Q)^{-1}Q^T F} + \\
& \mathbf{\mathrm{P}\{\mathbf{X}\} FQ (Q^T F^2Q)^{-1}Q^T } \Big ]_{ii}  
\end{split}
\label{eq:grad_F}
\end{equation}

\noindent
We use $\mathbf{F}^{-1}\mathrm{P}\{ \mathbf{FQ}\} = \mathbf{Q (Q^T F^2Q)^{-1} Q^T F} $ to ease the notation ( the diagonality of $\mathbf{F}$ guarantees it to be invertible). Our gradient function can be expressed as
\begin{equation}
    \begin{split}
\frac{\partial \nu(\mathbf{F}) }{\partial \mathbf{F}_{ii}} = &\frac{1}{3} \, \Big[
\mathbf{F}^{-1}\mathrm{P}\{ \mathbf{FQ}\} \, \mathrm{P}\{\mathbf{X}\} - 2 \mathbf{F}^{-1}\mathrm{P}\{ \mathbf{FQ}\} \,\mathrm{P}\{\mathbf{X}\} \, \mathrm{P}\{ \mathbf{FQ}\}+ \mathrm{P}\{\mathbf{X}\} \mathrm{P}\{ \mathbf{FQ}\} \mathbf{F}^{-1} \Big] _{ii}\\
    \end{split}
\end{equation}
Given $\mathbf{F}^{-1}\mathrm{P}\{ \mathbf{FQ}\} \, \mathrm{P}\{\mathbf{X}\}$ is symmetric, we have

\begin{equation}
\frac{\partial \nu(\mathbf{F}) }{\partial \mathbf{F}_{ii}} = \frac{2}{3} \, \Big[  \mathbf{F}^{-1}\mathrm{P}\{ \mathbf{FQ}\} \, \mathrm{P}\{\mathbf{X}\} - \mathbf{F}^{-1}\mathrm{P}\{ \mathbf{FQ}\} \,\mathrm{P}\{\mathbf{X}\} \, \mathrm{P}\{ \mathbf{FQ}\} \Big]_{ii}
\end{equation}

\noindent
This equation can be further merged into
\begin{equation}
\frac{\partial \nu(\mathbf{F}) }{\partial \mathbf{F}_{ii}} = \frac{2}{3} \, \Big[ \mathbf{F}^{-1}\mathrm{P}\{ \mathbf{FQ}\} \, \mathrm{P}\{\mathbf{X}\} (\mathbf{I} -\mathrm{P}\{ \mathbf{FQ}\})  \Big]_{ii}
\end{equation}
where $\mathbf{I}$ is the identity matrix.

Now, let us rewrite the derivative in terms of the underlying filter vector $\mathbf{f}$. First, remember that $\mathbf{F} = diag(\mathbf{f})$.  Let us denote the inverse operator - the one that extracts the diagonal from a square matrix and places the result in a vector - as $ediag$ (`e' signifies to `extract' the diagonal elements). Clearly, $ediag(diag(\mathbf{f}))=\mathbf{f}$.  Here, $diag()$ is a forward operation turning a vector into a diagonal matrix and $ediag$ is the companion reverse operator extracting the diagonal.

Now we can derive the gradient, in terms of the underlying filter vector, as
\begin{equation}
\nabla \nu ({\mathbf f}) =  \frac{\partial \nu(\mathbf{F}) }{\partial \mathbf{f}} = \frac{2}{3} \, ediag \Big( \mathbf{F}^{-1}\mathrm{P}\{ \mathbf{FQ}\} \, \mathrm{P}\{\mathbf{X}\} (\mathbf{I} -\mathrm{P}\{ \mathbf{FQ}\}) \Big)
\label{eq:grad_final}
\end{equation}
where the gradient with respect to the filter vector, $\nabla \nu ({\mathbf f})$, is a $n \times 1$  vector.  
\end{proof}

It is evident that the gradient has a very interesting structure: it is the diagonal of the product of three projection matrices multiplied by the inverse of the filter (at hand).  We will come back to this interesting feature later.

\subsection{Smoothness Constrained Filter}
When the filter is composed by a linear combination of a set of basis functions, $\mathbf{f} = \mathbf{Bc}$ where columns of $\mathbf{B}$ are basis vectors and $\mathbf{c}$ denotes the coefficients~\cite{EI2020}. By using the chain rule, we can calculate the gradient with respect to the coefficient vector $\mathbf{c}$ as:
\begin{equation}
\nabla \nu ({\mathbf c})     = \mathbf{B^T}  \frac{\partial \nu(\mathbf{F}) }{\partial \mathbf{f}}
\end{equation}
or, equivalently in its explicit form as
\begin{equation}
  \nabla \nu ({\mathbf c}) = \frac{2}{3} \,
  \mathbf{B^T}
  \, ediag \Big( \mathbf{F}^{-1}\mathrm{P}\{ \mathbf{FQ}\} \, \mathrm{P}\{\mathbf{X}\} (\mathbf{I} -\mathrm{P}\{ \mathbf{FQ}\}) \Big).
  \label{eq:gradient_Bc}
\end{equation}

\subsection{Filter Design with Regularization}
In Eq.~(\ref{eq:mu_regularization}), we reformulate the filter design optimization that we aim to minimize (instead of to maximize as for the Vora-Value optimization):
\begin{equation}
    \mu (\mathbf{F}) = - tr(\mathrm{P}\{\mathbf{ FQ}\}\mathrm{P}\{\mathbf{X}\}) + \;  \alpha ||\, {\mathbf F} \,|| ^2_2.
    \label{eq:mu_regularization}
\end{equation}
In contrast to the Vora-Value optimization, we reverse to the negative and cancel the fractional scalar in the equation. Here we use the symbol $\mu$ to denote the new objective function with a regularization term. The penalty term that we introduce here is the squared norm  of the filter where $\alpha > 0$.

Clearly,  we have the following relation of the gradient: $\nabla \mu ({\mathbf f}) =   -3 \nabla \nu ({\mathbf f}) + 2\alpha \mathbf{f}$. Therefore, the gradient of the regularized optimization is written as
\begin{equation}
\nabla \mu ({\mathbf f}) =
-2 \, ediag \Big( \mathbf{F}^{-1}\mathrm{P}\{ \mathbf{FQ}\} \, \mathrm{P}\{\mathbf{X}\} (\mathbf{I} -\mathrm{P}\{ \mathbf{FQ}\}) \Big) + 2\alpha \mathbf{f}.
\label{eq:grad_reg}
\end{equation}

\section{Derivation of Hessian Matrix}
Here we present how we derive  the second derivative - the Hessian matrix -  of our objective function given in Eq.~(\ref{eq:mu_regularization}).
The Hessian matrix makes a further derivative of Eq.~(\ref{eq:grad_reg}). As the second derivative will have many terms in the equation, to ease the notation, we use $\mathbf{A}$, $\mathbf{B}$ and $\mathbf{C}$ to respectively denote $\mathrm{P} \{ \mathbf{X} \}$, $\mathrm{P} \{ \mathbf{FQ} \}$ , $\mathbf{F^{-1}}$ hereafter. 

Now we differentiate for the second derivative using the matrix calculus laws in Eq.~(\ref{eq:calculus_rules}):
\begin{equation}
    \begin{split}
    \D ^2 \mu   = & \; tr(-2 \mathbf{  C B \,\D F\, CBA  \,\D F}  + 2 \mathbf{CB \,\D F\, CBAB \,\D F } \\
    & \mathbf{ + CBC \,\D F\, A  \,\D F  - CBC \,\D F\, AB  \,\D F }\\
    & \mathbf{- CBA \,\D F\, CB  \,\D F} + 2 \mathbf{CBAB \,\D F\, CB  \,\D F }\\
    & \mathbf{- CBABC  \,\D F  \,\D F } ) + 2\alpha  \, \mathbf{\D F \,\D F } 
    \end{split}
\end{equation}
After further merging between two terms in each line, we obtain
\begin{equation}
    \begin{split}
    \D ^2 \mu   = & tr( -2 \mathbf{ C B \,\D F\, CBA ( I- B)  \, \D F }\\
    & \mathbf{+ CBC  \,\D F  \, A(I - B)  \, \D F }\\
    & \mathbf{- CBA (I} - 2 \mathbf{B)  \,\D F  \, CB  \,\D F } \\
    & \mathbf{- CBABC  \, \D F  \, \D F } ) +  2\alpha \mathbf{ \D F \,\D F}
    \end{split}
    \label{eq:2nd_derivative_merge}
\end{equation}
where $\mathbf I$ denotes the $31\times 31$ identity matrix. 

Given $\mathbf{f} = diag(\mathbf{F})$ and any two square matrices $\mathbf{M}$ and $\mathbf{N}$, we have $tr(\mathbf{M}^T \,\D \mathbf{F\, N \,\D F} ) =  \sum_i \sum_j \mathbf{ M}_{i,j} \mathbf{N}_{j,i} \, \D \mathbf{f}_i \,\D \mathbf{f}_j$ where elements having the same indices in two matrices are multiplied. Using this property into Eq.~(\ref{eq:2nd_derivative_merge}) and the symmetric property of matrices $\mathbf{A,B,C}$ (as projector matrices and the  diagonal filter matrix $\mathbf{F}$  are symmetric), we can derive the Hessian matrix as
\begin{equation}
    \begin{split}
     \mathbf{ H} = \; &   - 2\mathbf{\Big(BC \Big) } \circ  \mathbf{\Big((I - B) ABC \Big) }  + \mathbf{\Big(CBC \Big)} \circ \mathbf{\Big((I-B)A \Big)} \\
     & + \mathbf{\Big((I - 2B)ABC \Big)} \circ \mathbf{\Big(BC \Big)} - \mathbf{\Big(CBABC \Big)} \circ  \mathbf{I}  + 2\alpha \mathbf{I}
\end{split}
\end{equation}
where $\circ$ denotes the Hadamard product (or elementwise product) of two matrices, i.e.\ $(\mathbf{M \circ N)}_{i,j} = \mathbf{M}_{i,j} \mathbf{N}_{i,j}$. 
The explicit expansion of the equation over projector matrices are written as

\begin{equation}
\begin{split}
\mathbf{H} = \;  & - 2\mathbf{\Big(\mathrm{P}\{ FQ\}F^{-1} \Big)} \circ \mathbf{\Big((I - \mathrm{P}\{ FQ\})\mathrm{P}\{ X\} \mathrm{P}\{ FQ\}  F^{-1} \Big) } \\
 & + \mathbf{\Big(F^{-1}\mathrm{P}\{ FQ\}F^{-1} \Big) } \circ \mathbf{\Big((I - \mathrm{P}\{ FQ\}) \mathrm{P}\{X\}\Big) } \\
   & + \mathbf{\Big((I - 2\mathrm{P}\{ FQ\})\mathrm{P}\{ X\}\mathrm{P}\{ FQ\}F^{-1} \Big)} \circ \mathbf{\Big(\mathrm{P}\{ FQ\}F^{-1} \Big) } \; \\
  & +\mathbf{\Big(F^{-1}\mathrm{P}\{ FQ\}\mathrm{P}\{ X\}\mathrm{P}\{ FQ\}F^{-1} \Big)} \circ  \mathbf{I} \\
  &+ 2\alpha \mathbf{I}
\end{split}
\label{eq:2nd_derivative}
\end{equation}
where the last term relates to the regularization term.

\section{Positive Definiteness of Hessian Matrix}
The gradient of the Vora-Value based objective function (when we discount the regularizer) has an interesting structure as given in Eq.~(\ref{eq:grad_final}). It is the product of three projector matrices and the inverse of the filter matrix.  We find that the filter vector, $\mathbf{f}$, and the gradient, $\nabla \nu(\mathbf{f})$, are perpendicular to each other. That is, when we multiply $\mathbf{f}^T$ to it, we have $\mathbf{f}^T \nabla \nu(\mathbf{f}) = \mathbf{f}^T ediag((\mathbf I - \mathrm P\{\mathbf{ FQ}\})\mathrm P\{\mathbf V\}\mathrm P\{ \mathbf{FQ}\}\mathbf F^{-1})  =0$. Therefore, given this property, we have
\begin{equation}
\begin{split}
\mathbf{f}^T  \nabla \mu ({\mathbf f}) & = 
\;\mathbf{f}^T( -3 \nabla \nu ({\mathbf f}) +  2\alpha \, \mathbf{f}) = 2\alpha \, \mathbf{f}^T \mathbf{f}
\end{split}
\label{eq:special_gradient}
\end{equation}

Now, if we make the derivative to the both sides of Eq.~(\ref{eq:special_gradient}) with respect to $\mathbf{f}$, we obtain
\begin{equation}
    \nabla \mu + (\nabla^2 \mu)\, \mathbf{f}  = 4 \alpha \, \mathbf{f}.
    \label{eq:special_hessian}
\end{equation}
If we multiply  $\mathbf{f}^T$ to this equation, we have $\mathbf{f}^T \nabla \mu + \mathbf{f}^T (\nabla^2 \mu)\, \mathbf{f}  = 4 \alpha \, \mathbf{f}^T \mathbf{f}$. From Eq.~(\ref{eq:special_gradient}), we know $\mathbf{f}^T \nabla \mu(\mathbf{f}) = 2\alpha \, \mathbf{f}^T \mathbf{f}$. Hence, we get 
\begin{equation}
    \mathbf{f}^T (\nabla^2 \mu)\, \mathbf{f}  = 2\alpha \, \mathbf{f}^T \mathbf{f} >0, \quad \text{if} \quad \alpha>0
\end{equation}
which  guarantees the Hessian to be positive definite under a positive regularizer $\alpha$ (and also a physically plausible filter is a non-zero vector, $\mathbf{f} >0$). The positive-definite  property ensures the Hessian matrix to be invertible and thus we can use the Newton's method - which involves the inverse of the Hessian matrix~\cite{convexBoyd} - for the Vora-Value based filter optimization.

\reftitle{References}


\begin{thebibliography}{999}

\bibitem{bibmatrix} J. R. Magnus and N. Heinz, Matrix differential calculus with applications in statistics and econometrics, John Wiley \& Sons, 2nd ed., 2007.

\bibitem{bibVV} P. L. Vora and H. J. Trussell, ``Measure of goodness of a set of color-scanning filters," \textit{JOSA A}, vol. 10, no. 7, pp. 1499-1508, 1993.

\bibitem{bibGolub} G. H. Golub and C. F.
V. Loan, Matrix computations, Johns Hopkins University Press, 4th ed., 2012.


\bibitem{bibLutherFilter} G. D. Finlayson, Y. Zhu, and H. Gong, ``Using a Simple Colour Pre-filter to Make Cameras More Colorimetric,'' in \textit{Color Imaging Conference}, no. 1, pp. 182-186, 2018.

\bibitem{EI2020} G. D. Finlayson, and Y. Zhu, ``An improved optimization method for finding a color filter to make a camera more colorimetric,” in \textit{Electronic Imaging: Color Imaging: Displaying, Processing, Hardcopy, and Applications proceedings}, 2020.

\bibitem{convexBoyd}
S. Boyd  and L. Vandenberghe, Convex Optimization, Cambridge university press, 2004.
\end{thebibliography}
\end{document}